\documentclass[12pt,reqno,a4paper]{amsart}
\usepackage{latexsym,amsmath,amsfonts,amssymb,amsthm}
\theoremstyle{plain}
\newtheorem{Thm}{Theorem}
\newtheorem{Lem}{Lemma}

\theoremstyle{definition}
\newtheorem*{Ack}{Acknowledgment}
\theoremstyle{remark}

\def\N{\mathbb N}

\def\S{\mathcal S}
\def\1{{\bf 1}}

\def\jacob #1#2{\genfrac{(}{)}{}{}{#1}{#2}}
\def\pmod #1{\ ({\rm{mod}}\ #1)}

\def\|{\ |\ }

\begin{document}
\title{Some congruences for trinomial coefficients}
\author{Hui-Qin Cao} \author{Hao Pan}
\begin{abstract} We prove several congruences for trinomial coefficients.\end{abstract}

\address{Department of Applied Mathematics, Nanjing Audit University, Nanjing 210029,. People's Public of China}
\email{caohq@nau.edu.cn}
\address{Department of Mathematics, Nanjing University, Nanjing 210093,
People's Republic of China}
\email{haopan79@yahoo.com.cn}
\maketitle
\section{Introduction}
\setcounter{Lem}{0}\setcounter{Thm}{0}\setcounter{Cor}{0}
\setcounter{equation}{0}

In \cite{PS}, Pan and Sun proved the following congruence on the sums of binomial coefficients:
\begin{equation}\label{psc}
\sum_{k=0}^{p-1}\binom{2k}{k+d}\equiv\jacob{p-d}{3}\pmod{p},
\end{equation}
where $p>3$ is a prime, $0\leqslant d\leqslant p-1$ and $\jacob{\ \;}{\ \;}$ is the Legendre symbol. They also proved that for prime $p>3$ and integer $0\leqslant d\leqslant p-1$,
\begin{equation}\label{pscc}
\sum_{k=1}^{p-1}\frac1k\binom{2k}{k+d}\equiv\begin{cases}d^{-1}(-1+2(-1)^d+3[3\mid p-d])\pmod{p},&\text{if }1\leqslant d\leqslant p,\\ 0\pmod{p},&\text{if }d=0,\end{cases}
\end{equation}
where $[A]=1$ or $0$ according to whether the assertion $A$ holds. Subsequently, Sun and Tauraso \cite{ST} extended (\ref{psc}) and showed that
\begin{equation}\label{stc}
\sum_{k=0}^{p-1}\binom{2k}{k+d}\equiv\jacob{p-d}{3}+2pS_d\pmod{p^2},
\end{equation}
where
$$
S_d=\sum_{0<k<d}\frac{(-1)^{k-1}}{k}\jacob{d-k}{3}.
$$

On the other hand, trinomial coefficients $\binom{n}{k}_2$ are given by
$$
(1+x+x^{-1})^n=\sum_{k=-n}^n\binom{n}{k}_2x^k.
$$
As G. E. Andrews mentioned \cite{A}, trinomial coefficients had been investigated by Euler. And Andrews and R. J. Baxter \cite{AB} found the $q$-analogues of trinomial coefficients play an important rule in the hard hexagon model. However, there is a similar congruence for trinomial coefficients:
\begin{equation}\label{sdcc}
\sum_{k=0}^{p-1}\binom{k}{d}_2\equiv\begin{cases}
(-1)^{\frac{p+d-1}{2}}\pmod{p}\quad&\text{if } d \text{ is even},\\
0\pmod{p}\quad&\text{if } d \text{ is odd},
\end{cases}
\end{equation}
where $p>3$ is prime and $0\leq d<p$. In fact, we shall prove the following stronger result.
\begin{Thm}\label{sdc2}
Let $p>3$ be a prime and let $d$ be an integer with $0\leq d\leq
p-1$. If $d$ is odd, then
\begin{align}\label{sddo}
\frac{1}{p}\sum_{k=0}^{p-1}\binom{k}{d}_2\equiv\frac{(-1)^{\frac{d+1}2}}{2}\bigg(\sum_{k=1}^{(d-1)/2}\frac{(-1)^k}{k}-3\sum_{\substack{1\leq k\leq (d-1)/2\\
3\mid k+p}}\frac{(-1)^k}{k}\bigg)\pmod{p}.
\end{align}
And if $d$ is even, then
\begin{align}\label{sdde}
&\frac{1}{p}\bigg((-1)^{\frac d2}\sum_{k=0}^{p-1}\binom{k}{d}_2-(-1)^{\frac{p-1}2}\bigg)\notag\\
\equiv&-3\jacob{-2}{p}\frac{\S_{(p-\jacob{3}p)/2}}{p}-2\sum_{0\leq j<d/2}\frac{(-1)^j}{2j+1}+3\sum_{\substack{0\leq j<d/2\\ 3\nmid p-2j-1}}\frac{(-1)^j}{2j+1}\pmod{p},
\end{align}
where the recurrence sequence $\{\S_n\}$ is defined as
$$
\S_0=0,\quad \S_1=1,\quad \S_{n+1}=4\S_n-\S_{n-1} \text{ for }n\geq 1.
$$
\end{Thm}
We also have a congruence for the alternating sums of trinomial coefficients.
\begin{Thm}\label{sda}
Let $p>3$ be a prime and let $d$ be an integer with $0\leq d\leq p-1$.
We have
\begin{equation}
\label{dd} \sum_{k=0}^{p-1}(-1)^k\binom{k}{d}_2\equiv
(-1)^{d-1}d+(-1)^dp\sum_{k=d+1}^{\lfloor (p+d-1)/2\rfloor}\frac
1k\binom{2k-d}{k+1}\pmod{p^2}.
\end{equation}
\end{Thm}
Unfornately, in general, it is not easy to compute
$$
\sum_{k=d+1}^{\lfloor (p+d-1)/2\rfloor}\frac
1k\binom{2k-d}{k+1}\pmod{p}.
$$
However, when $d=0,1,2$, we may get
\begin{Thm}\label{sdac}
For prime ${p>3}$ we have
\begin{equation}\label{d0}
\sum_{k=0}^{p-1}(-1)^k\binom{k}{0}_2\equiv p\frac{3\jacob
p3-1}{2}\pmod{p^2},
\end{equation}
\begin{equation}\label{d1}
\sum_{k=0}^{p-1}(-1)^k\binom{k}{1}_2\equiv 1-p\frac{3\jacob
p3-1}{2}\pmod{p^2}
\end{equation}
and
\begin{equation}\label{d2}
\sum_{k=0}^{p-1}(-1)^k\binom{k}{2}_2\equiv -2+3p\jacob p3\pmod{p^2}.
\end{equation}
\end{Thm}
Quiet recently, Z.-W. Sun told us that with help of the software
{\it Mathematica}, he found a similar congruence for the sum
$$
\sum_{k=0}^{p-1}\frac{T_k}{3^k}.
$$
Here we shall give a proof of Sun's congruence.
\begin{Thm}\label{sd3}
Let $p>3$ be a prime. Then
\begin{align*}
\sum_{k=0}^{p-1}\frac{T_k}{3^k}\equiv\begin{cases}p\pmod{p^2}\qquad\text{if } p\equiv1\pmod{3},\\
0\pmod{p^2}\qquad\text{if } p\equiv 2\pmod{3}.\end{cases}
\end{align*}
\end{Thm}

The proofs of the above theorems will be given in the next sections.

\section{Proof of Theorem \ref{sdc2}}
\setcounter{Lem}{0}\setcounter{Thm}{0}\setcounter{Cor}{0}
\setcounter{equation}{0}

\begin{Lem}\label{sd}
Suppose that $n>d\geq 0$. Then
\begin{equation}\label{sdi}
\sum_{k=0}^{n-1}\binom{k}{d}_2=\sum_{\substack{k=1\\ k+d\equiv
1\pmod{2}}}^{n}\binom{n}{k}\binom{k-1}{\frac{k+d-1}{2}}.
\end{equation}
\end{Lem}
\begin{proof} Let $[x^d]P(x)$ denote the coefficient of $x^d$ in the expansion of the polynomial $P(x)$. Since
$$\sum_{k=0}^{n-1}(1+x+x^{-1})^k=\frac{(1+x+x^{-1})^n-1}{x+x^{-1}}=\frac{1}{x^{n-1}}\cdot\frac{(1+x+x^2)^n-x^n}{1+x^2},$$
we have
\begin{align*}
\sum_{k=0}^{n-1}\binom{k}{d}_2=&\sum_{k=0}^{n-1}[x^d](1+x+x^{-1})^k
=[x^d]\frac{1}{x^{n-1}}\cdot\frac{(1+x+x^2)^n-x^n}{1+x^2}\\
=&[x^{n+d-1}]\frac{(1+x+x^2)^n-x^n}{1+x^2}=[x^{n+d-1}]\sum_{k=1}^n\binom{n}{k}(1+x^2)^{k-1}x^{n-k}\\
=&\sum_{k=1}^n\binom{n}{k}[x^{k+d-1}](1+x^2)^{k-1}=\sum_{\substack{k=1\\ k+d\equiv
1\pmod{2}}}^{n}\binom{n}{k}\binom{k-1}{\frac{k+d-1}{2}}.
\end{align*}
\end{proof}
Noting that
$$\binom{p}{k}\equiv 0\pmod{p}\qquad\text{and}\qquad\binom{p-1}{k}\equiv (-1)^k\pmod{p}$$
for $0<k<p$, (\ref{sdcc}) immediately follows from (\ref{sdi}).
\begin{Lem} Suppose that $p>3$ is prime and $0\leq m<(p-1)/2$. Then
$$
\sum_{j=1}^{(p-1)/2}\frac1j\binom{2j}{j+m}\equiv\frac{1}{m}(1-3[3\mid p+m])\pmod{p}.
$$
\end{Lem}
\begin{proof}
Let $$a_m=\sum_{j=1}^{(p-1)/2}\frac 1j\binom{2j}{j+m}.$$ 
Define the polynomials $\{v_n(x)\}_{n\in\N}$ by
$$
v_0(x)=2, v_1(x)=x, \text{ and  } v_{n+1}(x)=xv_n(x)-v_{n-1}(x),\quad
n=1,2,\ldots.
$$
Applying Theorem 2.1 in \cite{ST10}, we have
\begin{align*}
m\sum_{j=1}^{\frac{p-1}{2}}\frac
1j\binom{2j}{j+m}-v_m(-1)=&-\sum_{j=0}^{\frac{p-1}{2}+m}\binom{p+1}{j}v_{\frac{p+1}{2}+m-j}(-1)-2\binom{p}{\frac{p-1}{2}+m}\\
\equiv&-v_{\frac{p+1}{2}+m}(-1)-(p+1)v_{\frac{p-1}{2}+m}(-1)\\
\equiv&-v_{\frac{p+1}{2}+m}(-1)-v_{\frac{p-1}{2}+m}(-1)=v_{\frac{p+1}{2}+m+1}(-1)\pmod{p}.
\end{align*}
Since $v_n(-1)=3[3\mid n]-1$ for all $n\in\N$, we have
\begin{align*}
a_m=\sum_{j=1}^{\frac{p-1}{2}}\frac 1j\binom{2j}{j+m}\equiv&\frac 1m(v_m(-1)+v_{\frac{p+3}{2}+m}(-1))\\
=&\frac 1m(3[3\mid m]+3[3\mid \frac{p+3}{2}+m]-2)\\
=&\frac 1m(3[3\mid m]+3[3\|p-m]-2)\\
=&\frac 1m(1-3[3\mid p+m])\pmod{p}.
\end{align*}
\end{proof}
\begin{Lem}\label{j2} Suppose that $p>3$ is prime and $0\leq m<(p-1)/2$. Then
$$
\sum_{j=1}^{(p-3)/2}\frac{(-1)^j}{2j+1}\equiv\jacob{-1}p\frac{2^{p-1}-1}{2p}\pmod{p}.
$$
\end{Lem}
\begin{proof} Clearly,
$$\sum_{j=1}^{(p-3)/2}\frac{(-1)^j}{2j+1}\equiv\sum_{j=1}^{(p-3)/2}\frac{(-1)^j}{2j+1}\binom{p-1}{2i}=\frac1p\sum_{j=1}^{(p-3)/2}(-1)^j\binom{p}{2j+1}\pmod{p}.
$$
And
$$
\sum_{j=1}^{(p-3)/2}(-1)^j\binom{p}{2j+1}=(-1)^{\lfloor\frac{p-1}{4}\rfloor}2^{\frac{p-1}2}-(-1)^{\frac{p-1}2}=\jacob{-2}p\bigg(2^{\frac{p-1}2}-\jacob{2}p\bigg).
$$
Finally,
\begin{align*}
&2^{p-1}-1=\bigg(2^{\frac{p-1}2}-\jacob{2}p\bigg)\bigg(2^{\frac{p-1}2}+\jacob{2}p\bigg)\\
=&\bigg(2^{\frac{p-1}2}-\jacob{2}p\bigg)^2+2\jacob2p\bigg(2^{\frac{p-1}2}-\jacob{2}p\bigg)
\equiv2\jacob2p\bigg(2^{\frac{p-1}2}-\jacob{2}p\bigg)\pmod{p^2}.
\end{align*}
\end{proof}

\begin{Lem}\label{j23} Suppose that $p>3$ is prime. Then
$$
\sum_{\substack{0\leq j\leq (p-3)/2\\ 3\nmid p-2j-1}}\frac{(-1)^j}{2j+1}\equiv\jacob{-1}p\frac{2^{p-1}-1}{3p}-\jacob{-2}{p}\frac{\S_{(p-\jacob{3}p)/2}}{p}\pmod{p}.
$$
\end{Lem}
\begin{proof}
This is an immediate consequence of \cite[Corollary 3.3]{Sun}. 
\end{proof}
Now we are ready to prove (\ref{sddo}) and (\ref{sdde}).
\begin{proof}[Proof of (\ref{sddo})] Suppose that $d=2m+1$.
Let $$S_m=\frac{1}{p}\sum_{k=0}^{p-1}\binom{k}{2m+1}_2.$$ According
to Lemma \ref{sd},
\begin{align*}
S_m=\frac{1}{p}\sum_{\substack{k=1\\
2|k}}^{p}\binom{p}{k}\binom{k-1}{\frac{k+2m}{2}}=&\sum_{j=1}^{(p-1)/2}\frac 1{2j}\binom{p-1}{2j-1}\binom{2j-1}{j+m}\\
\equiv&-\frac 12\sum_{j=1}^{(p-1)/2}\frac
1j\binom{2j-1}{j+m}\pmod{p}.
\end{align*}
By (\ref{pscc}),
$$\sum_{j=1}^{p-1}\frac{1}{j}\binom{2j}{j}\equiv 0\pmod{p}.$$
As $$\binom{2j}{j}\equiv 0\pmod{p}$$ for $(p-1)/2<j<p$, we obtain that
\begin{align*}
S_0&\equiv-\frac 12\sum_{j=1}^{(p-1)/2}\frac
1j\binom{2j-1}{j}=-\frac 14\sum_{j=1}^{(p-1)/2}\frac
1j\binom{2j}{j}\\
&\equiv-\frac 14\sum_{j=1}^{p-1}\frac 1j\binom{2j}{j}\equiv
0\pmod{p}.
\end{align*}
For $1\leq m\leq(p-3)/2$,
\begin{align*}
S_m&\equiv-\frac 12\sum_{j=1}^{(p-1)/2}\frac
1j\bigg(\binom{2j}{j+m}-\binom{2j-1}{j+m-1}\bigg)\\
&=-\frac 12\sum_{j=1}^{(p-1)/2}\frac 1j\binom{2j}{j+m}+\frac
12\sum_{j=1}^{(p-1)/2}\frac 1j\binom{2j-1}{j+m-1}\\
&\equiv-\frac 12\sum_{j=1}^{(p-1)/2}\frac
1j\binom{2j}{j+m}-S_{m-1}\pmod{p}.
\end{align*}
Then
$$
S_m\equiv-S_{m-1}-\frac 12\sum_{j=1}^{(p-1)/2}\frac 1j\binom{2j}{j+m}\pmod{p},
$$
and therefore
\begin{align*}
&S_m\equiv(-1)^mS_0+\frac{(-1)^{m+1}}{2}\sum_{k=1}^m(-1)^k\sum_{j=1}^{(p-1)/2}\frac 1j\binom{2j}{j+k}\\
\equiv&\frac{(-1)^{m+1}}{2}\sum_{k=1}^m(-1)^k\sum_{j=1}^{(p-1)/2}\frac 1j\binom{2j}{j+k}\\
\equiv&\sum_{k=1}^m\frac{(-1)^k}{k}(1-3[3\mid p+k])=\sum_{k=1}^m\frac{(-1)^k}{k}-3\sum_{\substack{k=1\\
3\mid k+p}}^m\frac{(-1)^k}{k}\pmod{p}.
\end{align*}
This concludes the proof.\end{proof}

\begin{proof}[Proof of (\ref{sdde})]
Suppose that $d=2m$. Note that
$$
\binom{k+1}{2m+1}_2=\binom{k}{2m}_2+\binom{k}{2m+1}_2+\binom{k}{2m+2}_2.
$$
Hence
$$
\sum_{k=0}^{p-1}\binom{k}{2m}_2+\sum_{k=0}^{p-1}\binom{k}{2m+2}_2=\sum_{k=0}^{p-1}\binom{k+1}{2m+1}_2-\sum_{k=0}^{p-1}\binom{k}{2m+1}_2=\binom{p}{2m+1}_2.
$$
Let
$$
S_{m}=\sum_{k=0}^{p-1}\binom{k}{2m}_2.
$$
Then for every $j\geq 1$,
$$
S_m=(-1)^jS_{m+j}+\sum_{i=0}^{j-1}(-1)^i\binom{p}{2m+2i+1}_2
$$
In particular, when $j=(p-1)/2-m$, we have
$$
S_m=(-1)^{(p-1)/2-m}+\sum_{i=0}^{(p-3)/2-m}(-1)^i\binom{p}{2m+2i+1}_2.
$$
Now
\begin{align*}
&\frac1p\binom{p}{2m+2i+1}_2=\frac1p\sum_{k=0}^p(-1)^k\binom{p}{k}\binom{2(p-k)}{p-k-(2m+2i+1)}\\
=&\frac2{p-2m-2i-1}\binom{2p-1}{p-2m-2i-2}+\sum_{k=1}^{p-1}(-1)^k\frac1{k}\binom{p-1}{k-1}\binom{2(p-k)}{p-k-2m-2i-1}\\
\equiv&\frac{2}{2m+2i+1}+\sum_{k=1}^{p-1}\frac1{p-k}\binom{2(p-k)}{p-k-2m-2i-1}\pmod{p}.
\end{align*}
Theorefore
\begin{align*}
&\frac{S_m-(-1)^{(p-1)/2-m}}{p}\\
\equiv&\sum_{i=0}^{(p-3)/2-m}(-1)^i\bigg(\frac{2}{2m+2i+1}+\sum_{k=1}^{p-1}\frac1{k}\binom{2k}{k-2m-2i-1}\bigg)\\
\equiv&\sum_{i=0}^{(p-3)/2-m}(-1)^i\bigg(\frac{2}{2m+2i+1}+\sum_{k=1}^{p-1}\frac1{k}\binom{2k}{k+2m+2i+1}\bigg)\\
\equiv&\sum_{i=0}^{(p-3)/2-m}(-1)^i\bigg(\frac{2}{2m+2i+1}+\frac3{2m+2i+1}([3\nmid p-(2m+2i+1)]-1)\bigg)\\
=&2\sum_{i=0}^{(p-3)/2-m}\frac{(-1)^i}{2m+2i+1}-3\sum_{\substack{0\leq i\leq (p-3)/2-m\\ 3\nmid p-(2m+2i+1)}}\frac{(-1)^i}{2m+2i+1}\pmod{p}.
\end{align*}
By Lemma \ref{j2},
\begin{align*}
\sum_{i=0}^{(p-3)/2-m}\frac{(-1)^{m+i}}{2m+2i+1}\equiv\jacob{-1}p\frac{2^{p-1}-1}{2p}-\sum_{j=0}^{m-1}\frac{(-1)^j}{2j+1}\pmod{p}.
\end{align*}
And by Lemma \ref{j23},
\begin{align*}
&\sum_{\substack{0\leq i\leq (p-3)/2-m\\ 3\nmid p-(2m+2i+1)}}\frac{(-1)^{m+i}}{2m+2i+1}\\
\equiv&\jacob{-1}p\frac{2^{p-1}-1}{3p}-\jacob{-2}{p}\frac{\S_{(p-\jacob{3}p)/2}}{p}-\sum_{\substack{0\leq j\leq m-1\\ 3\nmid p+2j+1}}\frac{(-1)^j}{2j+1}\pmod{p}.
\end{align*}
We are done.
\end{proof}

\section{Proofs of Theorems \ref{sda} - \ref{sd3}}
\setcounter{Lem}{0}\setcounter{Thm}{0}\setcounter{Cor}{0}
\setcounter{equation}{0}
\begin{proof}[Proof of Theorem \ref{sda}]
Note that
$$
\sum_{k=0}^{p-1}(-1)^k\binom{k}{p-1}_2=(-1)^{p-1}\binom{p-1}{p-1}_2=1
$$
and
$$
\sum_{k=0}^{p-1}(-1)^k\binom{k}{p-2}_2=(-1)^{p-2}\binom{p-2}{p-2}_2+(-1)^{p-1}\binom{p-1}{p-2}_2=-1+p-1=p-2.
$$
So we may assume that $d<p-2$.
Since
\begin{align*}
\sum_{k=0}^{p-1}(-1)^k(1+x+x^{-1})^k=\frac{1+(1+x+x^{-1})^p}{1+(1+x+x^{-1})}=\frac{x^p+(1+x+x^2)^p}{x^{p-1}(1+x)^2},
\end{align*}
we have
\begin{align*}
\sum_{k=0}^{p-1}(-1)^k\binom{k}{d}_2=&\sum_{k=0}^{p-1}(-1)^k[x^d](1+x+x^{-1})^k=[x^d]\sum_{k=0}^{p-1}(-1)^k(1+x+x^{-1})^k\\
=&[x^{p+d-1}]\frac{x^p+(1+x+x^2)^p}{(1+x)^2}\\
=&[x^{d-1}]\frac
1{(1+x)^2}+[x^{p+d-1}]\sum_{k=0}^p\binom{p}{k}x^{2k}(1+x)^{p-k-2}\\
=&(-1)^{d-1}d+\sum_{k=d+1}^{\lfloor
(p+d-1)/2\rfloor}\binom{p}{k}\binom{p-k-2}{p+d-1-2k}\\
=&(-1)^{d-1}d+p\sum_{k=d+1}^{\lfloor (p+d-1)/2\rfloor}\frac
1k\binom{p-1}{k-1}\binom{p-k-2}{k-d-1}.
\end{align*}
Observe that $$\binom{p-1}{k-1}\equiv(-1)^{k-1}\pmod{p}$$ and
\begin{align*}
\binom{p-k-2}{k-d-1}&=\frac{(p-k-2)(p-k-3)\cdots(p-2k+d)}{(k-d-1)!}\\
&\equiv(-1)^{k-d-1}\frac{(k+2)(k+3)\cdots(2k-d)}{(k-d-1)!}\\
&=(-1)^{k-d-1}\binom{2k-d}{k-d-1}=(-1)^{k-d-1}\binom{2k-d}{k+1}\pmod{p}.
\end{align*}
So (\ref{dd}) is valid.
\end{proof}
\begin{proof}[Proof of Theorem \ref{sdac}] Let $C_n=\frac1{n+1}\binom{2n}{n}$ be the Catalan number.
Applying (\ref{dd}) with $d=0$, {d=1} and $d=2$ respectively, we get
\begin{align*}
\sum_{k=0}^{p-1}(-1)^k\binom{k}{0}_2\equiv
p\sum_{k=1}^{(p-1)/2}\frac
1k\binom{2k}{k+1}=p\sum_{k=1}^{(p-1)/2}C_k\pmod{p^2},
\end{align*}
\begin{align*}
\sum_{k=0}^{p-1}(-1)^k\binom{k}{1}_2\equiv 1-p\sum_{k=2}^{\lfloor
p/2\rfloor}\frac 1k\binom{2k-1}{k+1}=1-\frac
p2\sum_{k=2}^{(p-1)/2}\frac{k-1}{k}C_k\pmod{p^2}
\end{align*}
and
\begin{align*}
\sum_{k=0}^{p-1}(-1)^k\binom{k}{2}_2&\equiv
-2+p\sum_{k=3}^{(p+1)/2}\frac
1k\binom{2k-2}{k+1}=-2+p\sum_{k=3}^{(p+1)/2}C_{k-1}\frac{(k-2)(k-1)}{k(k+1)}\\
&=-2+p\bigg(\sum_{k=3}^{(p+1)/2}C_{k-1}-\sum_{k=3}^{(p+1)/2}C_{k-1}\frac{4k-2}{k(k+1)}\bigg)\\
&=-2+p\bigg(\sum_{k=2}^{(p-1)/2}C_k-\sum_{k=3}^{(p+1)/2}\frac{C_k}{k}\bigg)\pmod{p^2}.
\end{align*}
In \cite{PS}, Pan and Sun have proved that
\begin{align*}
\sum_{k=0}^{p-1}C_k\equiv\frac{3\jacob p3-1}{2}\pmod{p} \ \text{ and
}\  \sum_{k=1}^{p-1}\frac{C_k}k\equiv\frac32\big(1-\jacob
p3\big)\pmod{p}.
\end{align*}
Clearly,
\begin{align*}
C_0=C_1=1,\quad C_2=2,\quad C_{p-1}=\frac 1p\binom{2p-2}{p-1}\equiv
-1\pmod{p}
\end{align*}
and
\begin{align*}
C_k\equiv 0\pmod{p} \quad\text{for}\quad (p-1)/2<k<p-1.
\end{align*}
Therefore
\begin{align*}
\sum_{k=1}^{(p-1)/2}C_k\equiv\sum_{k=0}^{p-1}C_k\equiv\frac{3\jacob
p3-1}{2}\pmod{p}
\end{align*}
and
\begin{align*}
\sum_{k=1}^{(p-1)/2}\frac{C_k}k\equiv\frac{1-3\jacob p3}2\pmod{p}.
\end{align*}
Hence
\begin{align*}
\sum_{k=2}^{(p-1)/2}\frac{k-1}{k}C_k=\sum_{k=1}^{(p-1)/2}C_k-\sum_{k=1}^{(p-1)/2}\frac{C_k}{k}\equiv
3\jacob p3-1\pmod{p}
\end{align*}
and
\begin{align*}
\sum_{k=2}^{(p-1)/2}C_k-\sum_{k=3}^{(p+1)/2}\frac{C_k}{k}\equiv\frac{3\jacob
p3-1}{2}-1-\frac{1-3\jacob p3}2+2=3\jacob p3\pmod{p}.
\end{align*}
This yields (\ref{d0}), (\ref{d1}) and (\ref{d2}). We are done.
\end{proof}

\begin{proof}[Proof of Theorem \ref{sd3}]
Clearly
\begin{align*}
&\sum_{k=0}^{p-1}\bigg(\frac{(1+x+x^{-1}}{3}\bigg)^k=\frac{(1+x+x^{-1})^p/3^p-1}{(1+x+x^{-1})/3-1}\\
=&\frac{1}{3^{p-1}x^{p-1}}\cdot\frac{(1+x+x^2)^p-3^px^p}{(1-x)^2}\\
=&\frac{1}{3^{p-1}x^{p-1}}\cdot\frac{(1-x^3)^p-(3x(1-x))^p}{(1-x)^{p+2}}.
\end{align*}
 Then
\begin{align*}
&\sum_{k=0}^{p-1}\frac{T_k}{3^k}=[x^0]\sum_{k=0}^{p-1}\bigg(\frac{(1+x+x^{-1}}{3}\bigg)^k\\
=&\frac{1}{3^{p-1}}[x^{p-1}]\frac{(1-x^3)^p-(3x(1-x))^p}{(1-x)^{p+2}}\\
=&\frac{1}{3^{p-1}}[x^{p-1}]\frac{(1-x^3)^p}{(1-x)^{p+2}}\\
=&\frac{1}{3^{p-1}}\sum_{0\leq
k<p/3}\binom{p}{k}(-1)^k\binom{-(p+2)}{p-1-3k}(-1)^{p-1-3k}\\
=&\frac{1}{3^{p-1}}\sum_{0\leq
k<p/3}\binom{p}{k}(-1)^{p-1-3k}\binom{2p-3k}{p-1-3k}\\
=&\frac{1}{3^{p-1}}\binom{2p}{p-1}+\frac{1}{3^{p-1}}\sum_{0<
k<p/3}\binom{p}{k}(-1)^{k}\binom{2p-3k}{p-1-3k}\\
=&\frac{2p}{3^{p-1}(p-1)}\binom{2p-1}{p-2}+\frac{p}{3^{p-1}}\sum_{0<
k<p/3}\frac{(-1)^{k}}{k}\binom{p-1}{k-1}\binom{2p-3k}{p-1-3k}.
\end{align*}
It is known that $\binom{2p-1}{p-2}\equiv(-1)^{p-2}=-1\pmod{p}$ and
$$
\binom{2p-3k}{p-1-3k}=\binom{p+(p-3k)}{p-1-3k}\equiv\binom{p-3k}{p-1-3k}=p-3k\equiv-3k\pmod{p}.
$$
Therefore
\begin{align*}
\sum_{k=0}^{p-1}\frac{T_k}{3^k}\equiv&
2p+3p\lfloor\frac{p-1}{3}\rfloor\\
=&\begin{cases}p\pmod{p^2}\qquad\text{if } p\equiv 1\pmod{3},\\
0\pmod{p^2}\qquad\text{if } p\equiv 2\pmod{3}.\end{cases}
\end{align*}
We are done.
\end{proof}
\begin{Ack}
We are grateful to Professor Zhi-Wei Sun for his very helpful
suggestions on our paper.
\end{Ack}

\end{document}